\documentstyle [12pt,fleqn] {article}

\newcommand{\z}{\mbox{$\bar{z}$}}
\newcommand{\zz}{\mbox{$\bar{z_{0}}$}}
\newcommand{\psig}{\mbox{$\psi^{\gamma}(q)$}}
\newcommand{\psiz}{\mbox{$\psi^{\gamma}_z(q)$}}
\newcommand{\chig}{\mbox{$\chi^{\gamma}(q)$}}
\newcommand{\chiz}{\mbox{$\chi^{\gamma}_z(q)$}}
\newcommand{\rp}{\mbox{$\bf R^{+}$}}
\newtheorem{lemma}{Lemma}

\title {A calculation with a bi-orthogonal wavelet transformation
\thanks{\it This work was partially supported by
CONICET(Argentina).}}
\author{H.Falomir, M.A.Muschietti, E.M.Santangelo and J.Solomin\\
Facultad de
Ciencias Exactas, U.N.L.P.\\c.c.67, 1900 La Plata, Argentina.}
\date{ June 1, 1993}

\begin{document}

\maketitle
\begin{abstract}

We explore the use of bi-orthogonal basis for continuous wavelet
transformations, thus relaxing the so-called admissibility condition on
the analyzing wavelet. As an application, we determine the eigenvalues
and corresponding radial eigenfunctions of the Hamiltonian of relativistic
Hydrogen-like atoms.

{\em Pacs}: 02.30.+g, 03.65.-w, 03.65.Db
\end{abstract}

\newpage

\section {- Introduction}

 Wavelet transforms have been successfully used in Mathematics,
Physics
 and Engineering \cite{meyer,daube,chui}. In particular, in the context of
Quantum
 Mechanics,
 continuous wavelet transforms have proved very useful, giving
rise - for
 example - to entirely new approaches to problems with spherical
 symmetry.
 For nonrelativistic Hydrogen-like atoms, an adequate choice of
the analyzing
  wavelet reduces the radial Schr\"odinger equation to a first
order
 differential equation, and the analyticity of wavelet
 coefficients leads,
  in a straightforward manner, to the determination of the
eigenvalues and
  their corresponding eigenfunctions \cite{Paul}.

In this context, the selection of an analyzing wavelet is
constrained by the
  "admissibility condition", which guarantees the existence of an
inverse
  transform \cite{grossmann,Paul}.
  On the other hand, bi-orthogonal basis
have been
  introduced in the context of discrete \cite{cohen1,cohen2} as well as
continuous
  \cite{tchami} transforms.

However, in some cases, computational convenience may suggest
that the most
adequate "analyzing wavelet" be a non-admissible and even a
non-square-
integrable function. This is the case, for example, for
relativistic
Hydrogen-like atoms, as we will see later.

It is the aim of this paper  to extend the wavelet analysis
to such situations, where it is not possible to construct an
orthogonal
continuous basis of ${\bf L}^2$, via the "$ax+b$" transform of
the
analyzing wavelet. In order to get an invertible transformation,
we will
rather restrict ourselves to a subspace  of ${\bf L}^2$
(containing the
bounded eigenstates of the Hamiltonian to be treated) and
make use of bi-orthogonal continuous basis.

In section 2, we consider the space where this wavelet transform
is well
defined and some of its properties. We propose sufficient
conditions for
a function to belong to the space of wavelet coefficients. Such
conditions
are satisfied by the transformed eigenfunctions of the
relativistic
Hamiltonian to be treated later.

In section 3, the radial Dirac equation  is solved for
relativistic
Hydrogen-like atoms. As in the nonrelativistic case presented in
reference
\cite{Paul}, the analyticity of the space of coefficients is
shown to
determine the spectrum. Moreover, the aymptotic behaviour of
functions
in this space allows for a determination of the associated
eigenfunctions.

Finally, in section 4, we present some comments and conclusions.

\section {- The transformation}

Let us consider a function \psig, solution of
\begin{equation}
\left({{d}\over{dq}} + {{2-\gamma}\over{q}}\right) \psig =
-\psig\, ,
\end{equation}
with $q \in [0,\infty)$:
\begin{equation}
\psig = q^{\gamma - 2} e^{-q}\, .
\label{aw}
\end{equation}

For $\gamma>1$, \psig\  is an admissible wavelet \cite{grossmann}. So,
by
considering its $"ax + b"$ group transformation,
\begin{equation}
\psiz = a^{3/2} e^{ibq} \left[(aq)^{\gamma - 2} e^{-aq}\right]
,{\rm with} \ z = b + ia\ {\rm and} \ a>0\, ,
\end{equation}
a continuous orthogonal basis of ${\bf L}^2 (\rp, q^2 dq)$ can be
defined as
$\{\psiz\}$ (For definetness, we will consider the radial part of
three
dimensional problems).

Therefore, the wavelet coefficient of a function
$f(q) \in {\bf L}^{2}(\rp, q^{2} dq)$ is given by:
\begin{equation}
\left( \psiz , f(q) \right) = a^{\gamma - 1/2} F(\z)\, ,
\end{equation}
where
\begin{equation}
F(\z) = {\cal L}^{\gamma}\left(f(q)\right)(\z)
= \int_{0}^{\infty} dq\,  e^{-i \z q} q^{\gamma} f(q)
\label{efe}
\end{equation}
is an analytic function of the variable $\z$ in the lower
half-plane.
One then has the reconstruction formula:
\begin{equation}
f(q) = {{2^{2\gamma - 2}}\over{2 \pi \Gamma\left(2\gamma
-2\right)}}
\int_{\{ Im\,  z > 0 \}} d\mu_{L}(z)\,
\left( \psiz , f(q) \right) \psiz
\label{rec}
\end{equation}
with $d\mu_{L}(z)$ the left invariant measure of the $"ax + b"$
group:
\begin{equation}
d\mu_{L}(z) = {{da\,  db}\over{a^{2}}}\, .
\end{equation}
Moreover, the following equality holds:
\begin{equation}
\int_{0}^{\infty}dq\,  q^{2} |f(q)|^{2} =
{{2^{2 \gamma - 2}}\over{2 \pi \Gamma(2 \gamma - 2)}}
\int_{\{Im\,  z > 0\}}d\mu_{L}(z)\,  (Im\,  z)^{2 \gamma - 1}
|F(\z)|^{2} \, ,
\label{normas}
\end{equation}
which shows that $F(\z)$ belongs to a Bergman space ${\cal
B}_{2\gamma - 1}$
(see reference \cite{Paul}).

Now, for $1/2< \gamma <1$, the analyzing wavelet chosen is not an
admissible one \cite{grossmann}. So, in this range, it is not possible
to
construct an orthogonal basis leading to the reconstruction
formula
(\ref{rec}). Moreover, for $0< \gamma \leq 1/2$, \psig\  is not
even a square
integrable function, so that the integral in equation (\ref{efe})
doesn't
exist for an arbitrary $f(q) \in {\bf L}^2 (\rp, q^2 dq)$.

In what follows, we will be interested in showing that it is
still
possible to use the transform in equation (\ref{efe}) for solving
an
eigenvalue problem, provided certain regularity conditions are
satisfied
by its solutions. We will also analyze which properties of an
authentic
wavelet transform do still hold in such a situation.

To this end, we will introduce a bi-orthogonal continuous basis.
That is,
we will make use of different functions in the process of
analysis and
later reconstruction:
\begin{equation}
f(q) =
\int_{\{ Im\,  z > 0 \}} d\mu_{L}(z)\,
\left( \psiz , f(q) \right) \chiz \, ,
\end{equation}
where \chiz\ is obtained - through the action of the group "$ax +b$" - from
a function \chig satisfying:
\begin{equation}
\int_{0}^{\infty}dq\,  q\, \psig ^{*} \chig = {{1}\over{2\pi}}\, .
\end{equation}

Then, the following Lemmas hold:\bigskip
\begin{lemma}
Let $f(q) \in {\bf L}^{1}_{loc}(\rp, q^{\gamma} dq)\cap
{\bf L}^{2}\left((1,\infty), dq)\right)$,
with $0< \gamma <1$, and
consider $F(\z)$ as defined in equation(\ref{efe}). Then:

\bigskip\noindent
a) $F(\z)$ is an analytic function in the half-plane
${Im\,  \z < 0}$.
Moreover, \linebreak
$F(\z)
\begin{array}{c}
   \\
 \rightarrow \\
 ^{^{|Re\, z|\rightarrow \infty}}
\end{array}
 0$, with $Im\,  z =
a >0$, and $F(\z)
\begin{array}{c}
 \\
 \rightarrow \\
 ^{^{ Im\,  z \rightarrow \infty}}
 \end{array}
 0$.

\bigskip\noindent
b) If $f(q) \sim q^{\alpha - 1}$ ($\alpha \geq 0$) for $q \sim 0$
and
$f(q)$ is bounded when $q \rightarrow \infty$, then ${\cal
L}^{\gamma}$ transforms
the operator $q d/dq$ into the operator $-\z \partial / \partial
\z - (\gamma + 1)$.

\bigskip\noindent
c) If $f(q) \in {\bf L}^{2}(\rp, q^2 dq)$ then:
$\partial_{\z} F(\z) \in {\cal B}_{2(\gamma + 1) - 1}$, and
\[
\int_{0}^{\infty} dq\, q^{2} |f(q)|^{2} = \]
\begin{equation}
{{2^{2(\gamma + 1) - 2}}\over{2 \pi \Gamma\left( 2(\gamma + 1) -
2\right)}}
\int_{Im\, z > 0} d\mu_{L}(z)\, \left( Im\, z \right)^{2(\gamma
+ 1) - 1}
|\partial_{\z} F(\z)|^{2}\, .
\label{nor}
\end{equation}
\end{lemma}

\bigskip \noindent
{\bf Proof:}

\bigskip\noindent
a) The function
\begin{equation}
F(\z) = F(b - i a) = \int_{0}^{\infty} dq\, e^{- i b q}
q^{\gamma}
f(q) e^{- a q}
\label{F(b-ia)}
\end{equation}
is the Fourier transform of
$q^{\gamma} f(q) e^{-aq} \in {\bf L}^{1}(\rp,dq)$. So:
\begin{equation}
F(\z)
\begin{array}{c}
  \\
\rightarrow \\
 ^{{}^{|Re\, \z| \rightarrow \infty}}
\end{array}
0\, ,\ {\rm for}\  Im\, \z = - a <0\, .
\end{equation}
The analyticity of $F(\z)$ and the fact that $F(\z)\rightarrow 0$
for
$Im\, {\z} \rightarrow -\infty$ \linebreak
are direct consequences of its definition (see
equation(\ref{F(b-ia)})), since \linebreak $q^{\gamma} f(q)
e^{-aq} \in
{\bf L}^{1}(\rp, dq)$ for $a > 0$.

\bigskip\noindent b) Now,
\[
\int_{\varepsilon}^{\Lambda} dq\, e^{-i\z q} q^{\gamma}
\left[ q {{d}\over{dq}} f(q)\right] = \]\[
e^{-i\z q} q^{\gamma + 1} f(q) |_{\varepsilon}^{\Lambda} -
\int_{\varepsilon}^{\Lambda} dq\, {{d}\over{dq}} \left[e^{-i\z q}
q^{\gamma + 1}\right] f(q) \]
\begin{equation}
\begin{array}{c}
 \\ \rightarrow \\
{  ^{\Lambda \rightarrow \infty}_{\varepsilon \rightarrow 0}}
\end{array}
-\left[\z\partial_{\z} + \gamma +1 \right]
\int_{0}^{\infty} dq\, e^{-i\z q} q^{\gamma} f(q)\, ,
\end{equation}
since the integrated term vanishes under the assumption made on
the behavior
of $f(q)$, and $ e^{-i\z q} q^{\gamma + 1} f(q) \in {\bf
L}^{1}(\rp, dq)$.

\bigskip\noindent c) Notice that:
\begin{equation}
\partial_{\z}F(\z) = -i \int_{0}^{\infty} dq\, e^{-i\z q}
q^{\gamma + 1} f(q) = {\cal L}^{\gamma + 1}\big(f(q)\big)(\z)
\end{equation}
is the analytic factor of the wavelet coefficient of $f(q)$ with
respect to
the wavelet $\psi^{\gamma + 1}_{z}(q) \in {\bf L}^{2}(\rp, q^{2}
dq)$, wich
is admissible (since $\gamma + 1 > 1$). Then, from equation
(\ref{normas})
we inmediately get equation (\ref{nor}). $\Box$

\begin{lemma}
Let $F(\z)$ be an analytic function in the half-plane $\{Im\,\z<
0\}$,
with an asymptotic behaviour given by:
\begin{equation}
F(\z) = C_{0}\, (\z - \zz)^{-(\gamma + \alpha)} +
C_{1}\, (\z - \zz)^{-(\gamma + \alpha + 1)} +
G(\z)\, ,
\label{asym}
\end{equation}
where $ Im\, \zz >0$ and
$|G(\z)|\leq K\, |\z|^{-(\gamma + \alpha + 2)}$ is locally
bounded in the half-plane $Im \, \z\leq 0$ ($C_{0}$,
$C_{1}$ and $K$ are constants).
Then:

\bigskip\noindent
a) $(Im\, \z)^{\gamma - 1/2} F(\z)$ is the wavelet coefficient of
a function $f(q) \in$\linebreak
$ {\bf L}^{1}_{loc} (\rp,q^{\gamma} dq) \cap
{\bf L}^{2}\left((1,\infty), dq\right)$, given by:
\begin{equation}
f(q) = \int_{0}^{\infty} {{da}\over{a^2}}\,
\int_{-\infty}^{\infty} db\,\, (Im\, z)^{\gamma - 1/2} F(\z)
\chiz\, ,
\label{inverse}
\end{equation}
with $\z = b - ia$.

\bigskip\noindent
b) If $\partial_{\z} F(\z)\in {\cal B}_{2(\gamma + 1) - 1}$,
then $f(q)\in {\bf L}^{2}(\rp,q^{2} dq)$.

\bigskip\noindent
c) If
$|\z{{\partial}_{\z}}G(\z)|\leq K'\, |\z|^{-(\gamma + \alpha +
2)}$,
and is locally
bounded in the half-plane $Im \, \z\leq 0$ ($K'$ is a constant),
then
$\z{{\partial}_{\z}} F(\z) = {\cal L}^{\gamma}(h(q))(\z)$, where
\begin{equation}
h(q) = -(q\, {{d}\over{dq}} + \gamma + 1) f(q)\, .
\end{equation}
\end{lemma}

\noindent {\bf Proof:}

\bigskip\noindent
a) In the first place, notice that $(\z - \zz)^{-(\gamma +
\alpha)}$,
with $\alpha\geq 0$, is the analytic factor in the wavelet
coefficient
corresponding to the function
$f_{0}(q) = C_{0}
\left( {i^{(\gamma + \alpha)}/{\Gamma(\gamma + \alpha)}}  \right)
q^{\alpha - 1} e^{i\zz q}\in {\bf L}^{1}_{loc} (\rp,q^{\gamma}
dq)
\cap {\bf L}^{2}\left((1,\infty), dq\right)$. In fact,
\[
{\cal L}^{\gamma}[q^{\alpha -1} e^{i\zz q}](\z) =
{\cal F}[q^{\gamma + \alpha - 1} e^{-(a - a_{0})q} ](b - b_{0}) \]
\begin{equation}
 = \int_{0}^{\infty} dq\, q^{\gamma + \alpha - 1}
 e^{- i (\z - \zz)q} =
 {{\Gamma (\gamma + \alpha)}\over{[i (\z - \zz)]^{\gamma +
\alpha}}}\, .
\end{equation}
So:
\[ \int_{\{ Im\,  z > 0 \}} d\mu_{L}(z)\, (Im\,z)^{\gamma - 1/2} \,
(\z - \zz)^{-(\gamma +\alpha)} \chiz = \]
\begin{equation}
\int_{0}^{\infty} da\, a^{\gamma -1} \chi^{\gamma}(a q) 2 \pi
{\cal F}^{-1}[(\z - \zz)^{-(\gamma + \alpha)}](q) =
{{i^{\gamma + \alpha}}\over{\Gamma (\gamma + \alpha)}}
q^{\alpha - 1} e^{i \zz q}\, .
\label{reconstr}
\end{equation}
(Notice that the integral in the first member is conditionally
convergent).
A similar result holds, changing $\alpha$ into $\alpha + 1$, for
the second
term in equation (\ref{asym}), which is the analytic factor in the
wavelet
coefficient of $f_{1}(q) = C_{1} \left({i^{(\gamma + \alpha + 1)}
/{\Gamma(\gamma + \alpha + 1)}}\right) q^{\alpha} e^{i\zz q}
\in {\bf L}^{1}_{loc} (\rp,q^{\gamma} dq)\cap
{\bf L}^{2}\left((1,\infty), dq\right)$.

As regards $G(\z)$, under the assumptions made, it belongs to the
Bergman
space $B_{2\gamma + 1}$, since
\begin{equation}
\int_{\{ Im\,  z > 0 \}}  d\mu_{L}(z)\, (Im\,z)^{2(\gamma +1)-1}
|G(\z)|^{2} < \infty \, ,
\end{equation}
as can be easily verified: For example,
\begin{equation}
\int_{0}^{1} \int_{-1}^{1} |G(\z)|^{2} a^{2\gamma -1}\, db\,da <
\infty\, ,
\end{equation}
since $|G(\z)|$ is locally bounded.

Moreover, $\{ \psi_{z}^{\gamma + 1} \}$
is an orthogonal wavelet basis of ${\bf L}^{2}(\rp, q^{2} dq)$,
which defines
a bijection onto $B_{2\gamma + 1}$ (see reference \cite{Paul}). Then,
$g(q)\in
{\bf L}^{2}(\rp, q^{2} dq)$ exists such that:
\begin{equation}
G(\z) = \int_{0}^{\infty} dq\, q^{\gamma + 1} e^{- i \z q} g(q)
= {\cal L}^{\gamma +1}\big(g(q)\big)(\z) \, ,
\end{equation}
or, equivalently:
\begin{equation}
G(\z) = \int_{0}^{\infty} dq\, q^{\gamma} e^{- i \z q} f_{2}(q) =
{\cal F}[ q^{\gamma}  f_{2}(q)  e^{- a q}] =
{\cal L}^{\gamma}\big(f_{2}(q)\big)(\z) \, ,
\end{equation}
with $f_{2}(q) = q g(q) \in {\bf L}^{2}(\rp, dq)$ and, therefore,
$f_{2}(q)
\in {\bf L}^{1}_{loc} (\rp,q^{\gamma} dq)
\cap {\bf L}^{2}\left((1,\infty), dq\right)$.

Finally
\[
\int_{\{ Im\,  z > 0 \}} d\mu_{L}(z)\, (Im\,z)^{\gamma - 1/2} \,
G(\z) \chiz = \]
\[ \int_{0}^{\infty} da\, a^{\gamma -1} \chi^{\gamma}(a q)
\int_{-\infty}^{\infty} db\, G(b - i a) e^{i b q } = \]
\begin{equation}
2 \pi \int_{0}^{\infty} da\, a^{\gamma -1} \chi^{\gamma}(a q)
q^{\gamma} f_{2}(q) e^{- a q} = f_{2}(q) \, ,
\end{equation}
where use has been made of the fact that $G(\z)$ is the Fourier
transform
of a square integrable function.

\bigskip\noindent
b) Let us suppose that $\partial_{\z} F(\z) \in B_{2\gamma + 1}$;
then a function $h(q) \in {\bf L}^{2}(\rp,q^{2} dq)$ exists such
that:
\begin{equation}
 \partial_{\z} F(\z) =
 \int_{0}^{\infty} dq\, q^{\gamma + 1}  h(q) e^{- i \z q}
= {\cal L}^{\gamma +1}\big(h(q)\big)(\z) \, .
\end{equation}
Moreover, from a), we know that:
\[ \partial_{\z} F(\z) = \partial_{\z} \int_{0}^{\infty} dq\,
q^{\gamma}
[f_{0}(q) + f_{1}(q) + f_{2}(q) ] e^{- i \z q}   \]
\begin{equation}
 = -i  \int_{0}^{\infty} dq\, q^{\gamma +1}
 [f_{0}(q) + f_{1}(q) + f_{2}(q) ] e^{- i \z q} \, ,
\end{equation}
since the last integral is absolutely convergent.

Then, from a) (with $\gamma \rightarrow \gamma +1 $), we conclude that
$f_{0}(q) + f_{1}(q) + f_{2}(q) =$ \linebreak
$f(q) = i h(q) \in {\bf L}^{2}(\rp,q^{2} dq)$.

\bigskip\noindent
c) In the first place, we will consider, for $\alpha \geq 0$:
\[ \int_{\{ Im\,  z > 0 \}} d\mu_{L}(z)\, (Im\,z)^{\gamma - 1/2} \,
\left[ \z \partial_{\z} (\z - \zz)^{-(\gamma + \alpha)}\right]
\chiz =  \]
\begin{equation}
\hskip -1cm
-(\gamma + \alpha)
\int_{\{ Im\,  z > 0 \}} d\mu_{L}(z)\, (Im\,z)^{\gamma - 1/2} \,
\left[ {{1}\over{(\z - \zz)^{\gamma + \alpha}}} +
{{\zz}\over{(\z - \zz)^{\gamma + \alpha + 1} }}\right] \chiz \, .
\end{equation}
{}From equation (\ref{reconstr}), the previous expresion reduces to:
\[ -(\gamma + \alpha)  \left[{{i^{\gamma + \alpha}}\over{\Gamma
(\gamma + \alpha)}} q^{\alpha - 1} e^{ i \zz q} + \zz
{{i^{\gamma + \alpha + 1}}\over{\Gamma(\gamma + \alpha + 1)}}
q^{(\alpha + 1) - 1} e^{i \zz q}\right] =\]
\begin{equation}
- (q {{d}\over{dq}} + \gamma + 1) {{i^{\gamma +
\alpha}}\over{\Gamma
(\gamma + \alpha)}} q^{\alpha - 1} e^{ i \zz q}\, ,
\end{equation}
which proves the statement for the first two terms in equation
(\ref{asym}).

As concerns the third one:
\[
\int_{\{ Im\,  z > 0 \}} d\mu_{L}(z)\, (Im\,z)^{\gamma - 1/2} \,
[\z \partial_{\z} G(\z)] \chiz =\]
\[ \int_{0}^{\infty} da\, a^{\gamma -1} \chi^{\gamma}(a q)
\left[\int_{-\infty}^{\infty} db\,(b - i a) \partial_{b} G(b - i a)
e^{i b q }\right] = \]
\begin{equation}
\int_{0}^{\infty} da\, a^{\gamma -1} \chi^{\gamma}(a q)
\left[ -\int_{-\infty}^{\infty} db\, G(b - i a) (1 + a q + i b q)
e^{i b q }\right] \, ,
\label{third}
\end{equation}
where use has been made of the asymptotic behaviour of $G(\z)$ when
integrating by parts.

Notice that the integral between brackets in equation (\ref{third})
is
absolutely convergent, so that:
\[
\left[ \int_{-\infty}^{\infty} db\, G(b - i a) (1 + a q + i b q)
e^{i b q }\right] = \]
\begin{equation}
\left( 1 + a q + q {{d}\over{d q}} \right)
\int_{-\infty}^{\infty} db\, G(b - i a) e^{i b q}\, .
\end{equation}
Now, since $G(\z) \in B_{2\gamma + 1}$, one has:
\[ {\cal F}^{-1}\left[ G(b -i a)\right](q) = {{1}\over{2 \pi}}
\int_{-\infty}^{\infty} db\, G(b - i a) e^{i b q} \]
\begin{equation}
= q^{\gamma} g(q) e^{-a q}\, ,
\end{equation}
with $g(q) \in {\bf L}^{2}(\rp, dq)$. Therefore:
\[ \left( 1 + a q + q {{d}\over{d q}} \right)  2 \pi
q^{\gamma} g(q) e^{-a q} = \]
\begin{equation}
2 \pi q^{\gamma} e^{-a q}
\left( 1 + \gamma + q {{d}\over{d q}} \right) g(q)
\end{equation}
and
\[
\int_{\{ Im\,  z > 0 \}} d\mu_{L}(z)\, (Im\,z)^{\gamma - 1/2} \,
[\z \partial_{\z} G(\z)]\, \chiz =\]
\[ \int_{0}^{\infty} da\, a^{\gamma -1} \chi^{\gamma}(a q)
\left[ - 2 \pi q^{\gamma} e^{-a q}
\left( 1 + \gamma + q {{d}\over{d q}} \right) g(q) \right] \]
\begin{equation}
= - \left( 1 + \gamma + q {{d}\over{d q}} \right) g(q)  \, ,
\end{equation}
which completes the proof. $\Box$

\bigskip\noindent

For $0< \gamma <1$, the space of wavelet coefficients that appears
in Lemma 1.a) consists of
functions $(Im\, z)^{\gamma - 1/2} F(\z)$, where $F(\z)$ is
analytic for \linebreak
$Im\, \z<0$, vanishes for $\z\rightarrow \infty$ and is such that
$\partial_{\z} F(\z)$ belongs to
${\cal B}_{2(\gamma + 1)-1}$. This space of coefficients corresponds to the
transforms of functions in \linebreak
${\bf L}^{1}_{loc}(\rp, q^{\gamma} dq)\cap$
${\bf L}^{2}\left((1,\infty), dq\right)$.

\bigskip
Now, we introduce the linear space ${\cal A}_{\gamma}$ of functions
$F(\z)$,
analytic in the half plane
$Im\, \z<0$, vanishing for $\z\rightarrow \infty$ and such that
$\partial_{\z} F(\z) \in {\cal B}_{2(\gamma + 1)-1}$. Obviously, it
is a
pre-Hilbert space with respect to the scalar product:
\begin{equation}
<F|G>_{{\cal A}_{\gamma}} = \int_{Im z > 0} d\mu_{L}(z) (Im
z)^{2(\gamma+1)-1}
\partial_{\z}F(\z)^{*} \partial_{\z}G(\z)\, .
\end{equation}

\bigskip \begin{lemma}
The transformation ${\cal L}^{\gamma}$, defined in equation
(\ref{efe}) for
$0<\gamma<1$, maps a dense subspace of ${\bf L}^{2}(\rp, q^{2} dq)$
into
a dense subspace of the pre-Hilbert space
${\cal A}_{\gamma}$, preserving the norm.
\end{lemma}

\noindent {\bf Proof:}

Notice, in the first place, that the complete set of functions of
${\bf L}^{2}(\rp, q^{2} dq)$  given by
$\{\psi_{n}(q) = q^{\alpha - 1 + n}
 e^{-q}, n = 0, 1, 2...\}$, with $ 0\leq \alpha <1$,
is contained in
${\bf L}^{1}_{loc}(\rp, q^{\gamma} dq)\cap
{\bf L}^{2}\left((1,\infty), dq\right)$. So, ${\cal L}^{\gamma}$ is defined
on a dense subspace of ${\bf L}^{2}(\rp, q^{2} dq)$.
Moreover:
\[{\cal L}^{\gamma}\left( \psi_{n}(q) \right)(\z) =
\int_{0}^{\infty} dq\, q^{\gamma + \alpha -1 +n} e^{-q} e^{-i \z
q}\]
\begin{equation}
= \Gamma( \gamma + \alpha +n) \left[ i (\z - i)
\right]^{-(\gamma + \alpha +n)}\, .
\end{equation}

Now, the set
$\{ {\cal L}^{\gamma}\big(\psi_{n}(q)\big), n = 0, 1, 2...\}$ is
complete
in ${\cal A}_{\gamma}$, since
\begin{equation}
i \partial_{\z}{\cal L}^{\gamma}\left( \psi_{n}(q) \right)(\z) =
\int_{0}^{\infty} dq\, q^{\gamma + \alpha +n} e^{-q} e^{-i \z q} =
{\cal L}^{\gamma + 1}\left( \psi_{n}(q) \right)(\z) \, ,
\end{equation}
and because of the isometry established by the wavelet
transformation
${\cal L}^{\gamma + 1}$ between the Hilbert spaces
${\bf L}^{2}(\rp, q^{2} dq)$ and
${\cal B}_{2(\gamma + 1)-1}$ (see equation(\ref{normas})).

Finally, for $f(q), g(q) \in
{\bf L}^{2}(\rp, q^{2} dq) \cap
\big({\bf L}^{1}_{loc}(\rp, q^{\gamma} dq)\cap
{\bf L}^{2}\left((1,\infty), dq\right)\big)$ we have,
\[
<{\cal L}^{\gamma}\big(f(q)\big)(\z)|
 {\cal L}^{\gamma}\big(g(q)\big)(\z)>_{{\cal A}_{\gamma}} = \]
\[
<{\cal L}^{\gamma + 1}\big(f(q)\big)(\z)|
 {\cal L}^{\gamma + 1}\big(g(q)\big)(\z)>_{{\cal B}_{2\gamma + 1}} =
\]
\begin{equation}
 {{2 \pi \Gamma(2 \gamma - 2)}\over{2^{2 \gamma - 2}}}
\big( f, g\big)_{{\bf L}^{2}(\rp, q^{2} dq)}\, .\ \Box
\end{equation}

\bigskip\bigskip


\section{- Relativistic Hydrogen-like atom}

 As an application of the results presented in the previous
section,
we proceed, in what follows, to
the determination of the bounded eigenstates of the Hamiltonian of
relativistic Hydrogen-like atoms.

 As is well known \cite{Landau}, after elliminating angular
variables
through the SU(2) symmetry enjoyed by the problem at hand, the
radial
part of the eigenfunctions satysfies the following equations:
\[
{{df}\over{dr}} + {{1+\chi}\over{r}} f -
\left(\varepsilon + m + {{\lambda}\over{r}}\right) g = 0
\]
\begin{equation}
{{dg}\over{dr}} + {{1-\chi}\over{r}} g -
\left(\varepsilon - m + {{\lambda}\over{r}}\right) f = 0\, ,
\label{eqrad}
\end{equation}
where m is the electron mass, and $\varepsilon$ are the allowed
eigenvalues,
satisfying $|\varepsilon| < m$ for bounded states.

 Moreover, $\lambda = N\alpha$ (with $N$ the number of protons in
the
 nucleus and $\alpha = 1/137$, the fine structure constant). In
turn,
 $\chi$ is determined by the representation of SU(2) under study,
and
 is given by:
\begin{equation}
\chi = \left\{
\begin{array}{l}
+(j+1/2),\,  for\ j = l-1/2 \\
-(j+1/2),\  for\ j = l+1/2 \\
\end{array}  \right.\, ,
\end{equation}
with $j$ the total angular momentum of the electron.

 By defining:
\begin{equation}
0 \leq q = 2 r \sqrt{m^2 - \varepsilon^2},
\end{equation}
equation  (\ref{eqrad}) can be rewritten as:
\[
\left(q{{d}\over{dq}} + 1 + \chi\right) f(q)
- \left( {{d}\over{2}} \sqrt{{{m + \varepsilon}\over{m -
\varepsilon}}} +
\lambda \right) g(q) = 0
\]
\begin{equation}
\left(q{{d}\over{dq}} + 1 - \chi\right) g(q)
- \left( {{d}\over{2}} \sqrt{{{m + \varepsilon}\over{m -
\varepsilon}}} -
\lambda \right) f(q) = 0 \, ,
\label{ertrans}
\end{equation}
where $q\,  f(q)$ and $q\,  g(q)$ are square-integrable.

  As it can be easely seen \cite{Landau},
for $q\rightarrow 0$, the solutions of equation (\ref{ertrans})
behave as:
\begin{equation}
f(q),\ g(q) \sim  q^{-1 + \sqrt{\chi^2 - \lambda^2}}\, ,
\end{equation}
with $\chi^2 > \lambda^2$. So, $ f(q), g(q) \in
{\bf L}^{1}_{loc} (\rp,q^{\gamma} dq)$,
for $\gamma > 0$. The
transformation discussed in the previous section can therefore be
applied
since $f(q)$ and $ g(q)$ satisfy the requirements of
Lemma 1.

 Taking into account that the transformation is given by:
\begin{equation}
F(\z) = {\cal L}^{\gamma}\left(f(q)\right)(\z) =
\int_{0}^{\infty} dq \,  e^{-i \z q} q^{\gamma} f(q)\, ,
\end{equation}
it is easy to see that (Lemma 1):
\begin{equation}
{\cal L}^{\gamma}\,  q{{d}\over{dq}} = - \left(\z {{d}\over{d\z}}
+
\gamma +1 \right)\,  {\cal L}^{\gamma}\, ,
\end{equation}
and:
\begin{equation}
{\cal L}^{\gamma}\,  q = i{{d}\over{d\z}}\, {\cal L}^{\gamma}\, .
\end{equation}

 So, transforming equations (\ref{ertrans}), one gets:
\[
\left( - \z {{d}\over{d\z}} + \chi - \gamma\right) F(\z)
- \left( {{i}\over{2}} \sqrt{{{m + \varepsilon}\over{m -
\varepsilon}}}
\, {{d}\over{d\z}} + \lambda \right) G(\z) = 0
\]
\begin{equation}
\left( - \z {{d}\over{d\z}} - \chi - \gamma\right) G(\z)
- \left( {{i}\over{2}} \sqrt{{{m - \varepsilon}\over{m +
\varepsilon}}}
\, {{d}\over{d\z}} - \lambda \right) F(\z) = 0\, .
\label{eqz}
\end{equation}

 After some direct algebra, and calling
\begin{equation}
 \Phi(\z) =
 \left(\begin {array}{c} F(\z)\\ G(\z)\\
 \end{array}\right)\, ,
 \end{equation}
equation (\ref{eqz}) can be recast in the form:
\begin{equation}
{{d}\over{d\z}}\Phi(\z)  = -{{1}\over{2}}
\left\{
{{A' + B'}\over{\z - {{i}\over{2}} }} +
{{A' - B'}\over{\z + {{i}\over{2}}}}
\right\}
\Phi(\z)\, ,
\label{deriv}
\end{equation}
with:
\begin{equation}
A' = \left(
\begin{array}{cc}
\gamma - \chi & \lambda \\
-\lambda & \gamma + \chi\\
\end{array}
\right)
\end{equation}
\begin{equation}
B' = \left(
\begin{array}{cc}
\lambda
\sqrt{{{m + \varepsilon}\over{m - \varepsilon}}}  &
-(\chi + \gamma) \sqrt{{{m + \varepsilon}\over{m - \varepsilon}}}
\\  &  \\
-(\gamma - \chi)  \sqrt{{{m - \varepsilon}\over{m + \varepsilon}}}
&
- \lambda \sqrt{{{m - \varepsilon}\over{m + \varepsilon}}}   \\
\end{array}
\right)\, .
\end{equation}

 As is well known, the solution to equation (\ref{deriv}) is given
by:
\begin{equation}
\Phi(\z) = {\bf P} \exp \left\{
- {{1}\over{2}} \int_{\z_{0}}^{\z} d\z ' \,
\left[
{{A' + B'}\over{\z' - {{i}\over{2}} }} +
{{A' - B'}\over{\z' + {{i}\over{2}}}}
\right] \right\} \Phi(\z_{0})\, ,
\label{sol}
\end{equation}
where ${\bf P}$ means ordering over the path leading from $\z_0$
 to $\z$.

 Now, this expression can be greatly simplified through a judicious
choice of
 $\gamma$: By taking
 \footnote{Notice that, for $\chi^{2} = (j + 1/2)^{2} < 1 +
\lambda^{2}$,
 $\gamma < 1$, and we are in the conditions of the Lemmas of
Section 2.}
 \begin{equation}
 \gamma = + \sqrt {\chi^2 - \lambda^2} > 0\, ,
 \label{gamma}
 \end{equation}
 one has:
\[
(A')^{2} = 2 \gamma A' \ \  , \ \ A' B' =
{{2 \lambda\varepsilon}\over{\sqrt{m^2 - \varepsilon^2}}} A'\, ,\]
\begin{equation}
(B')^2 =
{{2 \lambda\varepsilon}\over{\sqrt{m^2 - \varepsilon^2}}} B' \ \ ,
\  \
B'A' = 2 \gamma B'\, ,
\end{equation}
and two new matrices can be defined as:
\begin{equation}
A = {{A' + B'}\over{- 2 \eta}} \ \ ,\ \
B = {{A' - B'}\over{- 2 \tilde{\eta}}}\, ,
\end{equation}
where
\begin{equation}
\eta = -\gamma - {{ \lambda\varepsilon}\over{\sqrt{m^2 -
\varepsilon^2}}}\, ,
\ \
\tilde{\eta} =
-\gamma + {{ \lambda\varepsilon}\over{\sqrt{m^2 -
\varepsilon^2}}}\, .
\end{equation}

 So, the following relations hold:
\[
A^2 = A \ \ ,\ \ AB = A\, ,\]
\begin{equation}
B^2 = B\ \ ,\ \ BA  = B \, .
\end{equation}

 For this choice of $\gamma$ it is easy to see that (\ref{sol})
reduces to:
\[ \Phi(\z) - \Phi(\z_0) =  \]
\begin{equation}
\int_{\z_0}^{\z} d\z' \,
\left(
{{\eta A}\over{\z' - {{i}\over{2}} }} +
{{\tilde{\eta} B}\over{\z' + {{i}\over{2}}}}
\right)
\left({{{\z' - {{i}\over{2}} } }\over{{\z_{0}' - {{i}\over{2}} } }}
\right)^{\eta}
\left({{{\z' + {{i}\over{2}} } }\over{{\z_{0}' + {{i}\over{2}} } }}
\right)^{\tilde \eta}
\Phi(\z_0)\, .
\end{equation}

\subsection*{Determination of the spectrum}

 As discussed in Section 2, $\Phi(\z)$ is an analytic function in
the lower
half-plane. So, its derivative:
\begin{equation}
{{d\Phi}\over{d\z}} =
\left(
{{\eta A}\over{\z' - {{i}\over{2}} }} +
{{\tilde{\eta} B}\over{\z' + {{i}\over{2}}}}
\right)
\left({{{\z' - {{i}\over{2}} } }\over{{\z_{0}' - {{i}\over{2}} } }}
\right)^{\eta}
\left({{{\z' + {{i}\over{2}} } }\over{{\z_{0}' + {{i}\over{2}} } }}
\right)^{\tilde \eta}
\Phi(\z_0)\, ,
\label{der2}
\end{equation}
must also be so. This requirement restricts $\tilde\eta$ to be a
nonnegative
integer:
\begin{equation}
\tilde\eta = - \gamma +
{{\lambda \varepsilon_{n}}\over{\sqrt{m^2-\varepsilon_{n}^2}}}
= n \, ,\ \ n = 0,1,...
\end{equation}
and $\eta = - n - 2\gamma$, from which the energy eigenvalues are
seen
to be:
\begin{equation}
{{\varepsilon_{n}}\over{m}} =
\left\{ 1 +
{{\lambda^2}\over{\left(\sqrt{\chi^2 - \lambda^2} + n \right)^2}}
\right\}^{-1/2} \, .
\end{equation}

 Thus, as in the nonrelativistic case \cite{Paul}, the bounded
spectrum
 can be determined from the requirement of analyticity on the
transform.

 \subsection*{Determination of eigenfunctions}

 From equation (\ref{der2}) and the condition $\Phi(\z) \rightarrow
0$ for
 $|\z| \rightarrow \infty$ (See Lemma 1 of Section 2) it can be
seen that:
\begin{equation}
\Phi(\z) \sim \z^{- 2 \gamma}\, , {\rm for}\  |\z| \rightarrow
\infty \, .
\end{equation}
 So, the limit:
\begin{equation}
\lim_{\z_{0} \rightarrow - i \infty}
{{\Phi(\z_0)}\over{
\left(\z_{0} - {{i}\over{2}}
\right)^{\eta}
\left(\z_{0} + {{i}\over{2}}
\right)^{\tilde \eta}
}} = \phi
\end{equation}
is finite.

 Moreover, for $\gamma$ as given in equation (\ref{gamma}), the
matrices $A$
 and $B$ can be written as:
\begin{eqnarray}
2 \eta A =
\left(-\gamma + \chi - \lambda
\sqrt{{{m + \varepsilon_n}\over{m - \varepsilon_n}}} \right)
\left(
\begin{array}{c}
1 \\ - \sqrt{{{m - \varepsilon_n}\over{m + \varepsilon_n}}}\\
\end{array}\right)
\otimes
\left( \begin{array}{cc}
1 & {{\lambda}\over{\gamma -\chi}}\\
\end{array} \right)\, ,  \\
2 \tilde{\eta}B = \left(-\gamma + \chi + \lambda
\sqrt{{{m + \varepsilon_n}\over{m - \varepsilon_n}}} \right)
\left(
\begin{array}{c}
1 \\ + \sqrt{{{m - \varepsilon_n}\over{m + \varepsilon_n}}}\\
\end{array}\right) \otimes
\left( \begin{array}{cc}
1 & {{\lambda}\over{\gamma -\chi}}\\
\end{array} \right)\, .
\end{eqnarray}

 Therefore, up to an overall multiplicative constant:
\[
\Phi_{n}(\z) =
\left(-\gamma + \chi - \lambda
\sqrt{{{m + \varepsilon_n}\over{m - \varepsilon_n}}} \right)
\left(
\begin{array}{c}
1 \\ - \sqrt{{{m - \varepsilon_n}\over{m + \varepsilon_n}}}\\
\end{array}\right)  \]\[
\int_{- i \infty}^{\z} d\z' \,
\left( \z' - {{i}\over{2}} \right)^{-(n+ 2 \gamma)-1}
\left( \z' + {{i}\over{2}} \right)^{n} + \]
\[
\left(-\gamma + \chi + \lambda
\sqrt{{{m + \varepsilon_n}\over{m - \varepsilon_n}}} \right)
\left(
\begin{array}{c}
1 \\ + \sqrt{{{m - \varepsilon_n}\over{m + \varepsilon_n}}}\\
\end{array}\right)  \]
\begin{equation}
\int_{- i \infty}^{\z} d\z' \,
\left( \z' - {{i}\over{2}} \right)^{-(n+ 2 \gamma)}
\left( \z' + {{i}\over{2}} \right)^{n-1} \, .
\label{phin}
\end{equation}

 The integrals in equations (\ref{phin}) can be evaluated on the
imaginary
 negative axis and analytically continued to the half-plane. In
this
 way, one obtains:
\[
\Phi_{n}(\z) =
\left(-\gamma + \chi - \lambda
\sqrt{{{m + \varepsilon_n}\over{m - \varepsilon_n}}} \right)
\left(
\begin{array}{c}
1 \\ - \sqrt{{{m - \varepsilon_n}\over{m + \varepsilon_n}}}\\
\end{array}\right)\]\[
\left( 1 + i\z \right)^{- 2 \gamma}
   \ _{2}F_{1}\left(-n,2 \gamma; 2\gamma +1;
{{1}\over{1+i\z}}\right) + \]
\[
\left(-\gamma + \chi + \lambda
\sqrt{{{m + \varepsilon_n}\over{m - \varepsilon_n}}} \right)
\left(
\begin{array}{c}
1 \\ + \sqrt{{{m - \varepsilon_n}\over{m + \varepsilon_n}}}\\
\end{array}\right) \]
\begin{equation}
\left(1 + i \z  \right)^{- 2 \gamma}
  \ _{2}F_{1}\left(-n+1,2 \gamma; 2\gamma +1;
{{1}\over{1+i\z}}\right) \, ,
\label{phin2}
\end{equation}
where $ _{2}F_{1}(...)$ is a Gauss hypergeometric function.

Notice that this solutions fulfill the hypothesis of Lemma 2, which
guarantees that we will obtain all the solutions in the
configuration space.
 In order to do so, the "inverse"
transform must be performed. To this end, an explicit $\chig$
function
must be chosen. For convenience, we adopt:
\begin{equation}
\chig = {{1}\over{2 \pi \Gamma (\gamma)}}\, .
\end{equation}
When inserting the first term in equation (\ref{phin2}) into equation
(\ref{inverse}), the integral
to be solved is then given by:
\[ \hskip -1cm
\lim_{N\rightarrow\infty} \int_{-N}^{N} db \,  \int_{0}^{\infty}
{{da}\over{a^2}}\,
a^{\gamma-1/2}\   {{a^{3/2} e^{i b q}}\over{2 \pi \Gamma (\gamma)}}
\left( 1 + i\z \right)^{- 2 \gamma}
   \ _{2}F_{1}\left(-n,2 \gamma; 2\gamma +1;
{{1}\over{1+i\z}}\right) \]
\[ \hskip -1cm
= {{1}\over{2 \pi \Gamma (\gamma)}}
\sum_{k=0}^{n} (-1)^{k}
{n \choose k}
{{2\gamma}\over{2\gamma+k}} \lim_{N\rightarrow\infty}
\int_{-N}^{N} db\,  e^{ibq} \int_{0}^{\infty} da\,  a^{\gamma -1}
\left[ 1+ib+a\right]^{-(2\gamma +k)}\]
\begin{equation}
= {{1}\over{\Gamma (2 \gamma)}} \theta(q) q^{\gamma -1} e^{-q}
\ _{1}F_{1}(-n, 2\gamma +1; q)\, ,
\label{hiper}
\end{equation}
where $ _{1}F_{1}(...)$ is a degenerate hypergeometric function.

 The second term in equation (\ref{phin2}) can similarly be
inverted
 (through the replacement $- n \rightarrow - n + 1$ in equation
(\ref{hiper})).
Thus, the eigenfunctions in the configuration space can be see to
coincide with the well known result (as given, for instance in reference
\cite{Landau}).

\section{- Conclusions}

In conclusion, we have explored the use of bi-orthogonal basis for
continuous wavelet transformations, a generalization which is aimed at
relaxing the so-called admissibility condition on the analyzing wavelet,
and turns out to be useful for computational reasons.

For definiteness, we have considered the radial dependence of functions in
${\bf R}^{3}$. As is well known, choosing as analyzing wavelet the function in
equation (\ref{aw}), with $\gamma > 1$, the wavelet transform in equation
(\ref{efe}) is an isometry between the Hilbert
spaces ${\bf L}^{2}(\rp, q^{2} dq)$ and ${\cal B}_{2\gamma - 1}$.

In Lemma 1, we have studied the transformation acting on functions
$f(q) \in {\bf L}^{1}_{loc}(\rp, q^{\gamma} dq)\cap
{\bf L}^{2}\left((1,\infty), dq)\right)$,
with $0< \gamma <1$, a region where the analyzing wavelet is not admissible
and can even be non square integrable. We have shown that the transform
$F(\z)$ so defined
is an analytic function in the half-plane
${Im\,  \z < 0}$,
such that
$F(\z)
 \rightarrow_{|Re\, z|\rightarrow \infty}
 0$, with $Im\,  z = a >0$, and $F(\z)
 \rightarrow_{ Im\,  z \rightarrow \infty} 0$,
and that the transformation maps differential operators acting on $f(q)$
into differential operators acting on $F(\z)$. Moreover, we have proved that,
if $f(q) \in {\bf L}^{2}(\rp, q^{2} dq)$, then $\partial_{\z} F(\z) \in
{\cal B}_{2\gamma + 1}$.

In Lemma 2, we have established that - for $F(\z)$ having an asymptotic
behaviour as given by equation (\ref{asym}) - the transformation has a
right inverse through the use of a bi-orthogonal basis.

In Lemma 3, we have shown that the transformation defined by equation
(\ref{efe}), for $0 < \gamma < 1$, is a mapping between a dense subspace of
${\bf L}^{2}(\rp, q^{2} dq)$  and a dense subspace of a pre-Hilbert space
${\cal A}_{\gamma}$, which preserves the norm (defined in ${\cal A}_{\gamma}$
in
terms of the
scalar product of derivatives in ${\cal B}_{2\gamma + 1}$).

Finally, as an example of the interest of our results, we have studied the
spectrum of relativistic Hydrogen-like atoms. We have shown that, in the
determination of eigenvalues of the Hamiltonian of this system and of their
associated radial eigenfunctions, a wavelet transformation can be employed,
and the calculation is greatly simplified by the choice
$\gamma = + \sqrt {\chi^2 - \lambda^2}$.
For physical reasons, $\gamma$ can be any real number greater than zero,
which makes apparent the need for our generalization of wavelet transforms.
By applying the results proved in our three Lemmas, we have determined the
spectrum from the requirement of analyticity on the transform, and we have
reconstructed the associated radial eigenfunctions through the use of a
bi-orthogonal basis. Both the eigenvalues and eigenfunctions thus obtained
can be seen to coincide with standard results.


\end{document}